\documentclass[10pt,a4paper,reqno]{amsart}
\usepackage{amsmath}
\usepackage{amsthm}
\usepackage{amsfonts}
\usepackage{amssymb}

\title{An $L^1$ estimate for half-space discrepancy}

\author{William W.L. Chen}

\address{Department of Mathematics, Macquarie University, Sydney, NSW 2109, Australia}

\email{wchen@maths.mq.edu.au}

\author{Giancarlo Travaglini}

\address{Dipartimento di Statistica, Edificio U7, Universit\`{a} di Milano-Bicocca, Via Bicocca degli Arcimboldi 8, 20126 Milano, Italy}

\email{giancarlo.travaglini@unimib.it}

\newtheorem{theorem}{Theorem}
\theoremstyle{remark}
\newtheorem*{remark}{Remark}
\theoremstyle{remark}
\newtheorem*{remarks}{Remarks}

\def\AA{\mathrm{A}}
\def\BB{\mathrm{B}}
\def\CC{\mathrm{C}}
\def\DD{\mathrm{D}}

\def\dd{\mathrm{d}}
\def\ee{\mathrm{e}}
\def\ii{\mathrm{i}}

\def\PPP{\mathcal{P}}
\def\TTT{\mathcal{T}}

\def\Rr{\mathbb{R}}
\def\Zz{\mathbb{Z}}

\def\card{\mathrm{card}}

\renewcommand{\le}{\leqslant}
\renewcommand{\ge}{\geqslant}

\begin{document}

\begin{abstract}
For every unit vector $\sigma\in\Sigma_{d-1}$ and every $r\ge0$, let
\begin{displaymath}
P_{\sigma,r}=[-1,1]^d\cap\{t\in\Rr^d:t\cdot\sigma\le r\}
\end{displaymath}
denote the intersection of the cube $[-1,1]^d$ with a half-space containing the origin $0\in\Rr^d$. We prove that if $N$ is the $d$-th power of an odd integer, then there exists a distribution $\PPP$ of $N$ points in $[-1,1]^d$ such that
\begin{displaymath}
\sup_{r\ge0}
\int_{\Sigma_{d-1}}\vert\card(\PPP\cap P_{\sigma,r})-N2^{-d}
\vert P_{\sigma,r}\vert\vert\,\dd\sigma
\le c_d(\log N)^d,
\end{displaymath}
generalizing an earlier result of Beck and the first author.
\end{abstract}

\maketitle

%
%

\section{Introduction}\label{sec:1}

The half-space discrepancy is a typical problem in the study of irregularities of point distribution, and represents a multi-dimensional variant of an open problem first posed by Roth; see Schmidt
\cite[pages 124--125]{schmidt77}. In its general form, it asks whether it is possible to choose $N$ points in a given bounded convex body in such a way that after cutting it into two parts by hyperplanes in different ways, the numbers of points in the two parts essentially depend only on the relative volumes. More precisely, let $\PPP$ denote a distribution of $N$ points in a bounded convex body $B\subset\Rr^d$. For every unit vector $\sigma\in\Sigma_{d-1}$ and every $r\ge0$, consider the half-space
$H_{\sigma,r}=\{t\in\Rr^d:t\cdot\sigma\le r\}$, where $\cdot$ denotes the usual inner product in $\Rr^d$, and let $S_{\sigma,r}=B\cap H_{\sigma,r}$. The problem is whether\footnote{We write $\vert S\vert$ to denote the Lebesgue measure of a Lebesgue measurable set $S$.}
\begin{equation}\label{eq:1.1}
\inf_{\card(\PPP)=N}\sup_{\substack{{r\ge0}\\{\sigma\in\Sigma_{d-1}}}}
\vert\card(\PPP\cap S_{\sigma,r})-N\vert B\vert^{-1}\vert S_{\sigma,r}\vert\vert
\end{equation}
is unbounded with $N$.

This question was first answered in the affirmative by Beck~\cite{beck83} in the case when $d=2$ and $B$ is the unit disc, using Fourier transform techniques. Subsequently, his almost sharp lower bound was improved by Alexander~\cite{alexander90} who used integral-geometric techniques to establish the $L^2$ result that for every distribution $\PPP$ of $N$ points in the unit disc, we have\footnote{Throughout this paper, the letter $c$ denotes positive absolute constants which may vary in value from one appearance to the next. Furthermore, the symbol $c$ with subscripts denotes positive constants whose values may depend on the subscripts displayed, and again may vary in value from one appearance to the next.}
\begin{equation}\label{eq:1.2}
\int_{\Sigma_1}\int_0^{\pi^{-1/2}}
\vert\card(\PPP\cap S_{\sigma,r})-N\vert S_{\sigma,r}\vert\vert^2\,\dd r\dd\sigma\ge cN^{1/2}.
\end{equation}
The unboundedness of \eqref{eq:1.1} in this special case follows immediately.

In fact, this last bound \eqref{eq:1.2} is sharp, in view of the amazing result of
Matou\v{s}ek~\cite{matousek95}, that there exist distributions $\PPP$ of $N$ points in the unit disc such that
\begin{displaymath}
\sup_{\substack{{r\ge0}\\{\sigma\in\Sigma_1}}}
\vert\card(\PPP\cap S_{\sigma,r})-N\vert S_{\sigma,r}\vert\vert
\le cN^{1/4},
\end{displaymath}
whereupon the upper bound
\begin{displaymath}
\int_{\Sigma_1}\int_0^{\pi^{-1/2}}
\vert\card(\PPP\cap S_{\sigma,r})-N\vert S_{\sigma,r}\vert\vert^2\,\dd r\dd\sigma\le cN^{1/2}
\end{displaymath}
follows immediately.

However, if one replaces the $L^2$ norm by the corresponding $L^1$ norm, one gets a rather different picture. No lower bound corresponding to \eqref{eq:1.2} is currently known, while Beck and the first author~\cite{BC93} have shown that for every bounded convex body $B\subset\Rr^2$ with centre of gravity at the origin and every natural number $N$, there exists a distribution $\PPP$ of $N$ points in
$B$ such that
\begin{equation}\label{eq:1.5}
\int_{\Sigma_1}\int_0^{R(\sigma)}
\vert\card(\PPP\cap S_{\sigma,r})-N\vert B\vert^{-1}\vert S_{\sigma,r}\vert\vert
\,\dd r\dd\sigma
\le c_B(\log N)^2,
\end{equation}
where $R(\sigma)=\sup\{t\cdot\sigma:t\in B\}$.

A careful description of the above and related problems can be found in Matou\v{s}ek
\cite[sections 3.2 and 6.6]{matousek99}.

The authors wish to express their gratitude to the referee for his careful reading of the manuscript and constructive comments.

%
%

\section{Main Results}\label{sec:2}

The purpose of this paper is to establish an estimate in the spirit of \eqref{eq:1.5}, in several variables and when the convex body is a cube. More precisely, let $Q=[-1,1]^d$. For every unit vector
$\sigma\in\Sigma_{d-1}$ and every $r\ge0$,
\begin{displaymath}
P_{\sigma,r}=Q\cap\{t\in\Rr^d:t\cdot\sigma\le r\}
\end{displaymath}
denotes the intersection of the cube $Q$ with one of the two half-spaces in $\Rr^d$ created by cutting
$\Rr^d$ by the hyperplane
\begin{equation}\label{eq:2.2}
s_{\sigma,r}=\{t\in\Rr^d:t\cdot\sigma=r\}.
\end{equation}
Our main result is the following.

\begin{theorem}\label{th:1}
Let $M>1$ be an integer, and let $N=(2M+1)^d$. Then there exists a distribution $\PPP$ of $N$ points in the cube $Q=[-1,1]^d$ such that
\begin{displaymath}
\sup_{r\ge0}\int_{\Sigma_{d-1}}
\vert\card(\PPP\cap P_{\sigma,r})-N2^{-d}\vert P_{\sigma,r}\vert\vert
\,\dd\sigma
\le c_d(\log N)^d.
\end{displaymath}
\end{theorem}

For every integer $M>1$, let
\begin{displaymath}
D_{\sigma,r}(M)=\card((M+\tfrac{1}{2})P_{\sigma,r}\cap\Zz^d)-(M+\tfrac{1}{2})^d\vert P_{\sigma,r}\vert.
\end{displaymath}
Theorem~\ref{th:1} follows immediately from the following result on lattice points by a simple scaling argument.

\begin{theorem}\label{th:2}
For every integer $M>1$, we have
\begin{displaymath}
\sup_{r\ge0}\int_{\Sigma_{d-1}}\vert D_{\sigma,r}(M)\vert\,\dd\sigma\le c_d(\log M)^d.
\end{displaymath}
\end{theorem}

\begin{remark}
If we replace the fraction $\frac{1}{2}$ in the dilation $(M+\frac{1}{2})P_{\sigma,r}$ of the set $P_{\sigma,r}$ by a different number in the interval $[0,1)$, then we obtain the trivial conclusion that $D_{\sigma,r}(M)$ is of order $M^{d-1}$ for every $\sigma$ and $r$.
\end{remark}

The remainder of the paper is organized as follows. In Section~\ref{sec:3}, we begin our proof of Theorem~\ref{th:2}, and split our argument into two cases. We then discuss these two cases separately in Sections \ref{sec:4} and~\ref{sec:5}.

%
%

\section{Fourier Transform and Divergence Theorem}\label{sec:3}

For every $x\in\Rr$, write
\begin{displaymath}
\varphi(x)
=(2-4\vert x\vert)_+
=\max\{2-4\vert x\vert,0\}.
\end{displaymath}
The function $\varphi$ is supported in $[-\frac{1}{2},\frac{1}{2}]$, satisfies
\begin{displaymath}
\int_\Rr\varphi(x)\,\dd x=1,
\end{displaymath}
and has Fourier transform $\widehat{\varphi}$ given by
\begin{equation}\label{eq:3.3}
\widehat{\varphi}(y)
=\int_{-1/2}^{1/2}\varphi(x)\,\ee^{-2\pi\ii yx}\,\dd x
=\left(\frac{2\sin(\pi y/2)}{\pi y}\right)^2.
\end{equation}
For every $t=(t_1,\ldots,t_d)\in\Rr^d$, write
\begin{displaymath}
\Phi(t)=\varphi(t_1)\ldots\varphi(t_d).
\end{displaymath}
For every positive integer $M$, write
\begin{displaymath}
\varphi_M(x)=M^{d-1}\varphi(M^{d-1}x)
\quad\mbox{and}\quad
\Phi_M(t)=M^{d^2-d}\Phi(M^{d-1}t).
\end{displaymath}
Then
\begin{displaymath}
\widehat{\Phi}_M(\xi)
=\int_{\Rr^d}\Phi_M(t)\,\ee^{-2\pi\ii\xi\cdot t}\,\dd t
=\widehat{\Phi}(M^{-d+1}\xi)
=\widehat{\varphi}(M^{-d+1}\xi_1)\ldots\widehat{\varphi}(M^{-d+1}\xi_d)
\end{displaymath}
for every $\xi=(\xi_1,\ldots,\xi_d)\in\Rr^d$. In particular, $\widehat{\Phi}_M(0)=1$.

Following a classical argument, we smooth the characteristic function
$\chi_{(M+\frac{1}{2})P_{\sigma,r}}$ of the set $(M+\frac{1}{2})P_{\sigma,r}$ by convolving it with
$\Phi_M$. We subsequently apply the Poisson summation formula to the convolution
\begin{displaymath}
\lambda_{M,(M+\frac{1}{2})P_{\sigma,r}}=\Phi_M\ast\chi_{(M+\frac{1}{2})P_{\sigma,r}}
\end{displaymath}
and deduce that
\begin{align}
\sum_{m\in\Zz^d}\lambda_{M,(M+\frac{1}{2})P_{\sigma,r}}(m)
&
=\sum_{m\in\Zz^d}\widehat{\lambda}_{M,(M+\frac{1}{2})P_{\sigma,r}}(m)
=\sum_{m\in\Zz^d}\widehat{\Phi}_M(m)\widehat{\chi}_{(M+\frac{1}{2})P_{\sigma,r}}(m)
\nonumber
\\
&
=(M+\tfrac{1}{2})^d
\sum_{m\in\Zz^d}\widehat{\Phi}_M(m)\widehat{\chi}_{P_{\sigma,r}}((M+\tfrac{1}{2})m).
\label{eq:3.8}
\end{align}

Observe next that the assumptions on $\Phi_M$ imply
\begin{align}
\sum_{m\in\Zz^d}\lambda_{M,(M+\frac{1}{2}-M^{-d+1})P_{\sigma,r}}(m)
&
\le\card((M+\tfrac{1}{2})P_{\sigma,r}\cap\Zz^d)
\nonumber
\\
&
\le\sum_{m\in\Zz^d}\lambda_{M,(M+\frac{1}{2}+M^{-d+1})P_{\sigma,r}}(m),
\nonumber
\end{align}
It follows from \eqref{eq:3.8} that
\begin{align}
D_{\sigma,r}(M)
&
\le\sum_{m\in\Zz^d}\lambda_{M,(M+\frac{1}{2}+M^{-d+1})P_{\sigma,r}}(m)
-(M+\tfrac{1}{2})^d\vert P_{\sigma,r}\vert
\nonumber
\\
&
=(M+\tfrac{1}{2}+M^{-d+1})^d
\sum_{m\in\Zz^d}\widehat{\Phi}_M(m)\widehat{\chi}_{P_{\sigma,r}}((M+\tfrac{1}{2}+M^{-d+1})m)
\nonumber
\\
&\qquad
-(M+\tfrac{1}{2})^d\vert P_{\sigma,r}\vert
\nonumber
\\
&
\le c_dM^d
\sum_{0\ne m\in\Zz^d}\widehat{\Phi}_M(m)\widehat{\chi}_{P_{\sigma,r}}((M+\tfrac{1}{2}+M^{-d+1})m)
+O(1),
\label{eq:3.10}
\end{align}
and a corresponding estimate holds from below.

We have to evaluate $\widehat{\chi}_{P_{\sigma,r}}(\xi)$ when $\vert\xi\vert\ge1$. By the divergence theorem, we have
\begin{displaymath}
\widehat{\chi}_{P_{\sigma,r}}(\xi)
=\int_{P_{\sigma,r}}\ee^{-2\pi\ii\xi\cdot t}\,\dd t
=\frac{\ii}{2\pi\vert\xi\vert^2}\int_{\partial P_{\sigma,r}}\ee^{-2\pi\ii\xi\cdot t}\,\xi\cdot\nu(t)\,\dd S_t,
\end{displaymath}
where $\nu(t)$ is the outward unit vector and $\dd S_t$ is the restriction of the Lebesgue measure to the boundary $\partial P_{\sigma,r}$, consisting of a bounded number of $(d-1)$-dimensional faces of $P_{\sigma,r}$. Let the polyhedron $G_{\sigma,r,d-1}$ denote one of these faces, and note that $\nu(t)$ is constant on $G_{\sigma,r,d-1}$. The study of $\widehat{\chi}_{P_{\sigma,r}}(\xi)$ therefore reduces to that of a finite number of terms of the form
\begin{equation}\label{eq:3.12}
\frac{\xi\cdot\nu}{\vert\xi\vert^2}\,\widehat{\mu}_{G_{\sigma,r,d-1}}(\xi),
\end{equation}
where $\mu_{G_{\sigma,r,d-1}}$ is the restriction of the Lebesgue measure to ${G_{\sigma,r,d-1}}$. We have two cases.

\textsc{Case $\AA_1$}. The face $G_{\sigma,r,d-1}$ is entirely contained in the hyperplane
$s_{\sigma,r}$; see \eqref{eq:2.2}.

\textsc{Case $\BB_1$}. The face $G_{\sigma,r,d-1}$ is entirely contained in one of the
$(d-1)$-dimensional faces of the cube $Q$.

%
%

\section{The Case $\AA_1$}\label{sec:4}

In this section, we consider the case when the face $G_{\sigma,r,d-1}$ is entirely contained in the hyperplane $s_{\sigma,r}$.

In this case, for every $r$, the face ${G_{\sigma,r,d-1}}$ rotates with
$\sigma\in\Sigma_{d-1}$, changing its shape as well as the number of its lower dimensional faces. However, the number of these lower dimensional faces and the lengths of their edges are bounded by positive constants that depend only on the dimension $d$. Then a mild variation of the proof of \cite[Theorem 2.1(ii)]{BCT97} gives
\begin{displaymath}
\int_{\Sigma_{d-1}}
\left\vert\frac{t\cdot\nu}{\vert t\vert^2}\,\widehat{\mu}_{G_{\sigma,r,d-1}}(t)\right\vert\,\dd\sigma
\le c_d\,\frac{(\log\vert t\vert)^{d-1}}{\vert t\vert^d},
\quad
\vert t\vert\ge2.
\end{displaymath}

\begin{remark}
According to \cite[Theorem 2.1(ii)]{BCT97}, a $d$-dimensional polyhedron $P$ satisfies
\begin{equation}\label{eq:4.2}
\int_{\Sigma_{d-1}}\vert\widehat{\chi}_P(\rho\sigma)\vert\,\dd\sigma
\le c_P\,\frac{(\log\rho)^{d-1}}{\rho^d},
\quad
\rho\ge2.
\end{equation}
If the diameter and the number of the faces of $P$ are bounded, then the constant $c_P$ can be replaced by a constant $c_d$. The proof of \eqref{eq:4.2} starts with the divergence theorem, and then proceeds by induction on the dimensions of the faces of $P$ in the following way. Write
$\sigma=(\cos\varphi,\eta\sin\varphi)\in\Sigma_{d-1}$, with $\eta\in\Sigma_{d-2}$ and $0\le\varphi\le\pi$. By the induction assumption, we have
\begin{displaymath}
\int_{\Sigma_{d-1}}\vert\widehat{\chi}_P(\rho\sigma)\vert\,\dd\sigma
\le c_P\,\frac{1}{\rho}\int_0^\pi\frac{(\log\rho\sin\varphi)^{d-2}}{(\rho\sin\varphi)^{d-1}}
\,(\sin\varphi)^{d-2}\,\dd\varphi
\le c_P\,\frac{(\log\rho)^{d-1}}{\rho^d},
\end{displaymath}
with the induction starting from the simple inequality
\begin{displaymath}
\int_0^\pi\frac{\vert\sin(\rho\sin\varphi)\vert}{\rho\sin\varphi}\,\dd\varphi
\le c\,\frac{\log\rho}{\rho}.
\end{displaymath}
\end{remark}

In the present case, the edges of $P_{\sigma,r}$ change in number and lengths under rotation and translation, but this does not affect the induction argument. For the first step, let $\gamma(\varphi,r)$ denote the length of a given edge on the boundary of $P_{\sigma,r}$. Then
\begin{displaymath}
\int_0^\pi\frac{\vert\sin(\rho\gamma(\varphi,r)\sin\varphi)\vert}{\sin\varphi}\,\dd\varphi
\le2+2\int_{1/\rho}^{\pi/2}\frac{1}{\sin\varphi}\,\dd\varphi
\le c\log\rho,
\end{displaymath}
and the contribution of Case $\AA_1$ to the estimate of
\begin{displaymath}
\int_{\Sigma_{d-1}}\vert D_{\sigma,r}(M)\vert\ \dd\sigma
\end{displaymath}
is bounded above (see \eqref{eq:3.10}) by
\begin{align}
&
c_dM^d\sum_{0\ne m\in\Zz^d}\vert\widehat{\Phi}_M(m)\vert
\,\frac{(\log M\vert m\vert)^{d-1}}{(M\vert m\vert)^d}
\nonumber
\\
&\quad
\le c_d(\log M)^{d-1}
\sum_{0\le\vert m\vert\le M^{d-1}}\vert\widehat{\Phi}_M(m)\vert\,\frac{1}{\vert m\vert^d}
+c_d\sum_{\vert m\vert>M^{d-1}}\vert\widehat{\Phi}_M(m)\vert
\,\frac{(\log\vert m\vert)^{d-1}}{\vert m\vert^d}
\nonumber
\\
&\quad
=\Theta_1+\Theta_2,
\label{eq:4.7}
\end{align}
say. Recall that the constants $c_d$ may change in value from one occurrence to the next.

We observe that $0\le\widehat{\Phi}_M(m)\le\widehat{\Phi}_M(0)=1$, and shall bound $\Theta_1$ by showing that
\begin{equation}\label{eq:4.8}
\sum_{\substack{{0<\vert m\vert\le M^{d-1}}\\{m_1\ge0,\ldots,m_d\ge0}}}\frac{1}{\vert m\vert^d}
\le c_d\log M,
\end{equation}
where $m=(m_1,\ldots,m_d)$. We shall achieve this by using induction to show that
\begin{equation}\label{eq:4.9}
\sum_{\substack{{0<\vert m\vert\le M^{d-1}}\\{m_1\ge0,\ldots,m_k\ge0}\\{m_{k+1}=0,\ldots,m_d=0}}}
\frac{1}{\vert m\vert^d}
\le c_d\log M
\end{equation}
holds for every $k=1,\ldots,d$.

Indeed, it is trivial to show that the inequality \eqref{eq:4.9} holds for $k=1$, noting that $d\ge2$. Suppose now that this inequality holds for every $k=1,2,\ldots,d-1$. Let
\begin{equation}\label{eq:4.10}
\TTT\overset{\text{def}}{=}
\left\{m=(m_1,\ldots,m_d)\in\Zz^d:m_1\ge1,\ldots,m_d\ge1,\max_jm_j\ge2\right\}.
\end{equation}
Then
\begin{equation}\label{eq:4.11}
\sum_{\substack{{0<\vert m\vert\le M^{d-1}}\\{m_1\ge0,\ldots,m_d\ge0}}}\frac{1}{\vert m\vert^d}
=\sum_{\substack{{0<\vert m\vert\le M^{d-1}}\\{\min_jm_j=0}}}\frac{1}{\vert m\vert^d}
+\frac{1}{d^{d/2}}
+\sum_{\substack{{0<\vert m\vert\le M^{d-1}}\\{m\in\TTT}}}\frac{1}{\vert m\vert^d}.
\end{equation}
Observe that the first sum on the right hand side of \eqref{eq:4.11} is a sum of a bounded number of terms of the form
\begin{equation}\label{eq:4.12}
\sum_{\substack{{0<\vert m\vert\le M^{d-1}}\\{m_{j_1}\ge0,\ldots,m_{j_k}\ge0}\\
{m_j=0\text{ if }j\not\in\{j_1,\ldots,j_k\}}}}
\frac{1}{\vert m\vert^d}
\end{equation}
with $k<d$. Their overall contribution does not exceed $c_d\log M$ by the induction hypothesis -- note that the quantity \eqref{eq:4.12} is invariant under permutation of the variables $m_1,\ldots,m_d$, and is therefore equal to the left hand side of \eqref{eq:4.9}. To study the last term on the right hand side of \eqref{eq:4.11}, we consider the bijection (see
\eqref{eq:4.10})
\begin{displaymath}
\TTT\ni m=(m_1,\ldots,m_d)
\longleftrightarrow
(m_1-1,m_1]\times\ldots\times(m_d-1,m_d]
\overset{\text{def}}{=}Q_m,
\end{displaymath}
and note that the union of the cubes $Q_m$ satisfies
\begin{displaymath}
\bigcup_{m\in\TTT}Q_m=(0,+\infty)^d\setminus(0,1]^d.
\end{displaymath}
Then
\begin{align}
\sum_{\substack{{0<\vert m\vert\le M^{d-1}}\\{m\in\TTT}}}\frac{1}{\vert m\vert^d}
&
\le\sum_{\substack{{0<\vert m\vert\le M^{d-1}}\\{m\in\TTT}}}\int_{Q_m}\frac{1}{\vert x\vert^d}\,\dd x
\le\int_{1\le\vert x\vert\le M^d}\frac{1}{\vert x\vert^d}\,\dd x
\nonumber
\\
&
=c_d\int_1^{M^d}\frac{1}{s}\,\dd s
=c_d\log M.
\nonumber
\end{align}
This completes the proof of the inequality \eqref{eq:4.8}.

We now conclude from \eqref{eq:4.7} and \eqref{eq:4.8} that
\begin{equation}\label{eq:4.16}
\Theta_1\le c_d(\log M)^d.
\end{equation}

To study the term $\Theta_2$, note first of all that for every $m=(m_1,\ldots,m_d)$, there exists at least one index $j^\ast$ such that $\vert m_{j^\ast}\vert\ge\vert m\vert/\sqrt{d}$. In view of \eqref{eq:3.3}, we have
\begin{align}
\widehat{\Phi}_M(m)
&
=\prod_{j=1}^d\widehat{\varphi}(M^{-d+1}m_j)
=\prod_{j=1}^d\left(\frac{2\sin(\pi M^{-d+1}m_j/2)}{\pi M^{-d+1}m_j}\right)^2
\nonumber
\\
&
\le\left(\frac{2\sin(\pi M^{-d+1}m_{j^\ast}/2)}{\pi M^{-d+1}m_{j^\ast}}\right)^2
\le c_dM^{2d-2}\frac{1}{\vert m\vert^2}.
\nonumber
\end{align}
It follows that
\begin{equation}\label{eq:4.18}
\Theta_2\le c_dM^{2d-2}\sum_{\vert m\vert>M^{d-1}}\frac{(\log\vert m\vert)^{d-1}}{\vert m\vert^{d+2}}.
\end{equation}
For $s>M^{d-1}$, the function
\begin{displaymath}
s\mapsto\frac{(\log s)^{d-1}}{s^{d+2}}
\end{displaymath}
decreases with $s$. We can then apply the earlier argument and control the right hand side of
\eqref{eq:4.18} with an integral, which can then be handled using integration by parts $d-1$ times. More precisely, we have
\begin{align}
\Theta_2
&
\le c_dM^{2d-2}\int_{M^{d-1}}^{+\infty}\frac{(\log s)^{d-1}}{s^3}\,\dd s
\nonumber
\\
&
\le c_dM^{2d-2}\left(M^{2-2d}(\log M)^{d-1}+\int_{M^{d-1}}^{+\infty}\frac{(\log s)^{d-2}}{s^3}\,\dd s\right)
\nonumber
\\
&
\le c_dM^{2d-2}\left(M^{2-2d}(\log M)^{d-1}+M^{2-2d}(\log M)^{d-2}
+\int_{M^{d-1}}^{+\infty}\frac{(\log s)^{d-3}}{s^3}\,\dd s\right)
\nonumber
\\
&
\le\ldots
\le c_d(\log M)^{d-1}.
\label{eq:4.20}
\end{align}

Combining \eqref{eq:4.7}, \eqref{eq:4.16} and \eqref{eq:4.20}, we conclude that the contribution of Case $\AA_1$ to the estimate of
\begin{displaymath}
\int_{\Sigma_{d-1}}\vert D_{\sigma,r}(M)\vert\ \dd\sigma
\end{displaymath}
is bounded above by $c_d(\log M)^d$.

%
%

\section{The Case $\BB_1$}\label{sec:5}

In this section, we consider the case when the face $G_{\sigma,r,d-1}$ is entirely contained in one of the
$(d-1)$-dimensional faces of the cube $Q$. Our proof is inductive in nature.

Without loss of generality, we may assume that $\nu=(0,\ldots,0,1)$, so that the face $G_{\sigma,r,d-1}$ is contained in the hyperplane $t_d=1$. Then \eqref{eq:3.12} becomes
\begin{align}
\frac{\xi\cdot\nu}{\vert\xi\vert^2}\,\widehat{\mu}_{G_{\sigma,r,d-1}}(\xi)
&
=\frac{\xi_d}{\vert\xi\vert^2}\,\ee^{-2\pi\ii\xi_d}
\int_{F_{\sigma,r,d-1}}\ee^{-2\pi\ii(\xi_1,\ldots,\xi_{d-1})\cdot(t_1,\ldots,t_{d-1})}\,\dd t_1\ldots\dd t_{d-1}
\nonumber
\\
&
=\frac{\xi_d}{\vert\xi\vert^2}\,\ee^{-2\pi\ii\xi_d}
\widehat{\chi}_{F_{\sigma,r,d-1}}(\xi_1,\ldots,\xi_{d-1}),
\label{eq:5.1}
\end{align}
where $F_{\sigma,r,d-1}=G_{\sigma,r,d-1}-\nu$ can be interpreted as a polyhedron in $\Rr^{d-1}$, with characteristic function $\chi_{F_{\sigma,r,d-1}}$. To study \eqref{eq:5.1}, we consider two cases.

\textsc{Case $\CC_1$}. We have $\vert(\xi_1,\ldots,\xi_{d-1})\vert<1$.

\textsc{Case $\DD_1$}. We have $\vert(\xi_1,\ldots,\xi_{d-1})\vert\ge1$.

We begin with Case $\CC_1$. Recall that $m\in\Zz^d$, so $\vert(m_1,\ldots,m_{d-1})\vert<1$ clearly implies $m_1=\ldots=m_{d-1}=0$. The contribution of this case to an upper estimate for \eqref{eq:3.10} therefore does not exceed
\begin{align}
&
M^d\left\vert\sum_{0\ne m_d\in\Zz}\widehat{\varphi}\left(\frac{m_d}{M^{d-1}}\right)\frac{1}{Mm_d}
\,\ee^{-2\pi\ii(M+\frac{1}{2}+M^{-d+1})m_d}\,\vert F_{\sigma,r,d-1}\vert\right\vert
\nonumber
\\
&\quad
=2\vert F_{\sigma,r,d-1}\vert\left\vert\sum_{m_d=1}^{+\infty}(-1)^{m_d}
\widehat{\varphi}\left(\frac{m_d}{M^{d-1}}\right)
\frac{\sin2\pi M^{-d+1}m_d}{M^{-d+1}m_d}\right\vert
\le H+K,
\label{eq:5.2}
\end{align}
where
\begin{displaymath}
H=2^d\left\vert\sum_{m_d=1}^{M^{d-1}}(-1)^{m_d}
\widehat{\varphi}\left(\frac{m_d}{M^{d-1}}\right)
\frac{\sin2\pi M^{-d+1}m_d}{M^{-d+1}m_d}\right\vert
\end{displaymath}
and
\begin{displaymath}
K=2^d\left\vert\sum_{m_d=M^{d-1}+1}^{+\infty}(-1)^{m_d}
\widehat{\varphi}\left(\frac{m_d}{M^{d-1}}\right)
\frac{\sin2\pi M^{-d+1}m_d}{M^{-d+1}m_d}\right\vert.
\end{displaymath}

\begin{remark}
Note that the equality in \eqref{eq:5.2} depends on the fraction $\frac{1}{2}$ in the dilation
$(M+\frac{1}{2})P_{\sigma,r}$ of the set $P_{\sigma,r}$.
\end{remark}

For the sum $H$, note that we have $0\le M^{-d+1}m_d\le1$, and that we can split the interval
$0\le x\le1$ into a bounded number of subintervals where the function
\begin{displaymath}
x\longrightarrow\widehat{\varphi}(x)\,\frac{\sin2\pi x}{x}
\end{displaymath}
is monotone and does not change sign. It follows that the sum $H$ is not greater than the sum of a bounded number of Leibniz sums, and this implies $H\le c_d$.

For the sum in $K$, note that
\begin{align}
&\sum_{m_d=M^{d-1}+1}^{+\infty}(-1)^{m_d}
\widehat{\varphi}\left(\frac{m_d}{M^{d-1}}\right)
\frac{\sin2\pi M^{-d+1}m_d}{M^{-d+1}m_d}
\nonumber
\\
&\qquad
=\sum_{j=1}^{+\infty}\sum_{m_d=jM^{d-1}+1}^{(j+1)M^{d-1}}(-1)^{m_d}
\widehat{\varphi}\left(\frac{m_d}{M^{d-1}}\right)
\frac{\sin2\pi M^{-d+1}m_d}{M^{-d+1}m_d}
\nonumber
\\
&\qquad
=\sum_{j=1}^{+\infty}(-1)^{jM^{d-1}}\sum_{m_d=1}^{M^{d-1}}(-1)^{m_d}
\widehat{\varphi}\left(j+\frac{m_d}{M^{d-1}}\right)
\frac{\sin2\pi M^{-d+1}m_d}{j+M^{-d+1}m_d}.
\nonumber
\end{align}
It is an exercise in the calculus to show the existence of a positive constant $c$ such that for every index~$j$, there are at most $c$ subintervals of the interval $0\le x\le1$ where the function
\begin{displaymath}
x\longrightarrow\widehat{\varphi}(j+x)\,\frac{\sin2\pi x}{j+x}
\end{displaymath}
is monotone and does not change sign. Then
\begin{displaymath}
K
\le c_d\sum_{j=1}^{+\infty}j^{-3}
\le c_d.
\end{displaymath}

Next, we turn our attention to Case $\DD_1$. Applying the divergence theorem to the polyhedron $F_{\sigma,r,d-1}$, we meet cases similar to Case $\AA_1$ and Case $\BB_1$. At the $\ell$-th step, where $1\le\ell\le d-1$, the divergence theorem leads to one of the following two cases.

\textsc{Case $\AA_\ell$}. We have a $(d-\ell)$-dimensional face entirely contained in some hyperplane in $\Rr^{d-\ell+1}$ analogous to $s_{\sigma,r}$.

\textsc{Case $\BB_\ell$}. We have a $(d-\ell)$-dimensional face entirely contained in one of the
$(d-\ell)$-dimensional faces of the cube $[-1,1]^{d-\ell+1}$.

In Case $\AA_\ell$, we proceed as in Case $\AA_1$.

In Case $\BB_\ell$, we need to study terms of the form
\begin{equation}\label{eq:5.9}
M^d\sum_{0\ne m\in\Zz^d}\left(\prod_{j=0}^{\ell-1}\Xi_j(m)\right)
\widehat{\chi}_{F_{\sigma,r,d-\ell}}((M+\tfrac{1}{2}+M^{-d+1})(m_1,\ldots,m_{d-\ell})),
\end{equation}
where, for $0\le j\le\ell-1$,
\begin{displaymath}
\Xi_j(m)=\widehat{\varphi}\left(\frac{m_{d-j}}{M^{d-1}}\right)
\!\frac{(M+\frac{1}{2}+M^{-d+1})m_{d-j}}{\vert(M+\frac{1}{2}+M^{-d+1})(m_1,\ldots,m_{d-j})\vert^2}
\,\ee^{-2\pi\ii(M+\frac{1}{2}+M^{-d+1})m_{d-j}}.
\end{displaymath}
We split the sum \eqref{eq:5.9} into the following two cases.

\textsc{Case $\CC_\ell$}. We have $\vert(m_1,\ldots,m_{d-\ell})\vert<1$.

\textsc{Case $\DD_\ell$}. We have $\vert(m_1,\ldots,m_{d-\ell})\vert\ge1$.

Since $m_1=\ldots=m_{d-\ell}=0$ in Case $\CC_\ell$, the contribution of this case to the sum \eqref{eq:5.9} is equal to
\begin{displaymath}
A
=M^d\vert F_{\sigma,r,d-\ell}\vert
\sum_{0\ne m_{d-\ell+1}\in\Zz}\ldots\sum_{0\ne m_d\in\Zz}
\left(\prod_{j=0}^{\ell-1}\Xi_j(m)\right).
\end{displaymath}
To study this sum, write
\begin{align}
&
B(m_{d-\ell+1},\ldots,m_{d-1})
=\sum_{0\ne m_d\in\Zz}\Xi_0(m)
\nonumber
\\
&
=\frac{2M^{-d+1}}{M+\frac{1}{2}+M^{-d+1}}\!\sum_{m_d=1}^{+\infty}
\widehat{\varphi}\left(\frac{m_d}{M^{d-1}}\right)\!(-1)^{m_d}
\frac{M^{-d+1}m_d\sin2\pi M^{-d+1}m_d}{\vert M^{-d+1}(0,\ldots,0,m_{d-\ell+1},\ldots,m_d)\vert^2},
\nonumber
\end{align}
and observe that the function
\begin{displaymath}
x\longrightarrow\widehat{\varphi}(x)
\,\frac{x\sin2\pi x}{\vert(0,\ldots,0,M^{-d+1}m_{d-\ell+1},\ldots,M^{-d+1}m_{d-1},x)\vert^2}
\end{displaymath}
is bounded in $x$, uniformly for $m_{d-\ell+1},\ldots,m_{d-1}$. Applying the earlier argument for the sum $H$ to each of these functions, we conclude that
\begin{displaymath}
\vert M^dB(m_{d-\ell+1},\ldots,m_{d-1})\vert\le c_d.
\end{displaymath}
Then
\begin{displaymath}
A
=M^d\vert F_{\sigma,r,d-\ell}\vert
\sum_{0\ne m_{d-\ell+1}\in\Zz}\ldots\sum_{0\ne m_{d-1}\in\Zz}
\left(\prod_{j=1}^{\ell-1}\Xi_j(m)\right)B(m_{d-\ell+1},\ldots,m_{d-1})
\end{displaymath}
satisfies
\begin{align}
\vert A\vert
&
\le c_d\sum_{0\ne m_{d-\ell+1}\in\Zz}\ldots\sum_{0\ne m_{d-1}\in\Zz}
\left(\prod_{j=1}^{\ell-1}\vert\Xi_j(m)\vert\right)
\nonumber
\\
&
\le c_d\sum_{0\ne m_{d-\ell+1}\in\Zz}\ldots\sum_{0\ne m_{d-1}\in\Zz}
\left(\prod_{j=1}^{\ell-1}
\left\vert\widehat{\varphi}\left(\frac{m_{d-j}}{M^{d-1}}\right)\right\vert
\frac{1}{Mm_{d-j}}\right)
\nonumber
\\
&
\le c_dM^{-\ell+1}
\left(\sum_{k=1}^{+\infty}
\left\vert\widehat{\varphi}\left(\frac{k}{M^{d-1}}\right)\right\vert\frac{1}{k}\right)^{\ell-1}
\nonumber
\\
&
\le c_dM^{-\ell+1}
\left(\sum_{k=1}^{M^{d-1}}\frac{1}{k}+\sum_{k=M^{d-1}+1}^{+\infty}\frac{M^{2d-2}}{k^3}\right)^{\ell-1}
\nonumber
\\
&
\le c_dM^{-\ell+1}(\log M)^\ell,
\label{eq:5.16}
\end{align}
using \eqref{eq:3.3}.

\begin{remarks}
(1) The upper bound \eqref{eq:5.16} is more than we need. However, the problem of bounding the sum $A$ is not entirely trivial, since simply putting absolute values inside the sums
$B(m_{d-\ell+1},\ldots,m_{d-1})$ does not lead to a useful estimate. Again we have used the cancellations given by the term $\frac{1}{2}$ in the dilation $(M+\frac{1}{2})P_{\sigma,r}$ of the set $P_{\sigma,r}$.

(2) It appears that we are studying the contribution of the boundary of $Q$ to the discrepancy, but we know that this contribution is actually zero, as a consequence of the term $\frac{1}{2}$ in the dilation, at least as far as whole faces of $Q$ are concerned. The delicate point here is that we are not estimating the actual discrepancy arising from the boundary, but have arrived at the boundary through the Poisson summation formula and the divergence theorem. Thus this approach does not seem to allow us to state mathematically that the contribution of the boundary must be negligible.
\end{remarks}

In case $\DD_\ell$, we again apply the divergence theorem, and meet cases similar to Case $\AA_\ell$ and Case $\BB_\ell$. At the last step, we have part of an edge of $Q$, say
$\{t,1,\ldots,1\}_{b(\sigma,r)\le t\le1}$. Then we need to bound the sum
\begin{equation}\label{eq:5.17}
M^d\sum_{0\ne m\in\Zz^d}\widehat{\Phi}_M(m)
\left(\prod_{j=2}^d\Upsilon_j(m)\right)
\int_{b(\sigma,r)}^1\ee^{-2\pi\ii(M+\frac{1}{2}+M^{-d+1})m_1s}\,\dd s,
\end{equation}
where, for $2\le j\le d$,
\begin{displaymath}
\Upsilon_j(m)
=\frac{(M+\frac{1}{2}+M^{-d+1})m_j}{\vert(M+\frac{1}{2}+M^{-d+1})(m_1,\ldots,m_j)\vert^2}
\,\ee^{-2\pi\ii(M+\frac{1}{2}+M^{-d+1})m_j}.
\end{displaymath}
The part of the sum \eqref{eq:5.17} where $m_1=0$ is Case $\CC_{d-1}$. For $m_1\ne0$, we compute the integral and bound the sum \eqref{eq:5.17}, uniformly in $\sigma$ and $t$, by
\begin{align}
M^d\prod_{j=1}^d
\left(\sum_{m_j=1}^{+\infty}\left\vert\widehat{\varphi}\left(\frac{m_j}{M^{d-1}}\right)\right\vert
\frac{1}{Mm_j}\right)
&
=\left(\sum_{k=1}^{+\infty}
\left\vert\widehat{\varphi}\left(\frac{k}{M^{d-1}}\right)\right\vert\frac{1}{k}\right)^\ell
\nonumber
\\
&
\le c_dM^{-d}(\log M)^d,
\nonumber
\end{align}
as in \eqref{eq:5.16}.

This completes the proof of Theorem~\ref{th:2}.

%
%

\end{document}